\documentclass[12pt]{amsart}

\usepackage[latin1]{inputenc}
\usepackage[T1]{fontenc}
\usepackage{hyperref}
\usepackage{geometry} 
\usepackage{stmaryrd}
\usepackage{amsmath,amssymb,amsthm,latexsym}

\newtheorem{theorem}{Theorem}[section]
\newtheorem{proposition}[theorem]{Proposition}
\newtheorem{fact}[theorem]{Fact}
\newtheorem{cor}[theorem]{Corollary}
\newtheorem{lmm}[theorem]{Lemma}
\newtheorem*{claim}{Claim}
\newtheorem{clm}{Claim}

\theoremstyle{definition}
\newtheorem{definition}[theorem]{Definition}
\newtheorem{remark}[theorem]{Remark}
\newtheorem{conj}[theorem]{Conjecture}
\newtheorem{expl}[theorem]{Example}

\def\bsp{\begin{expl}}
\def\ebsp{\end{expl}}
\def\behe{\begin{clm}}
\def\ebehe{\end{clm}}
\def\beh{\begin{claim}}
\def\ebeh{\end{claim}}
\def\defn{\begin{definition}}
\def\edefn{\end{definition}}
\def\satz{\begin{theorem}}
\def\esatz{\end{theorem}}
\def\tats{\begin{fact}}
\def\etats{\end{fact}}
\def\kor{\begin{cor}}
\def\ekor{\end{cor}}
\def\bema{\begin{remark}}
\def\ebema{\end{remark}}
\def\lem{\begin{lmm}}
\def\elem{\end{lmm}}
\def\bem{\begin{remark}}
\def\ebem{\end{remark}}
\def\verm{\begin{conj}}
\def\everm{\end{conj}}
\def\bew{\begin{proof}}
\def\ebew{\end{proof}}
\def\bewbeh{\begin{proof}[Proof of Claim]}
\def\satzli{\begin{proposition}}
\def\esatzli{\end{proposition}}

\def\Ind#1#2{#1\setbox0=\hbox{$#1x$}\kern\wd0\hbox to 
0pt{\hss$#1\mid$\hss}
\lower.9\ht0\hbox to 0pt{\hss$#1\smile$\hss}\kern\wd0}

\def\Notind#1#2{#1\setbox0=\hbox{$#1x$}\kern\wd0\hbox to 0pt{\mathchardef
\nn="3236\hss$#1\nn$\kern1.4\wd0\hss}\hbox to 0pt{\hss$#1\mid$\hss}\lower.9\ht0
\hbox to 0pt{\hss$#1\smile$\hss}\kern\wd0}

\def\indd{\mathop{\ \hbox to 
0pt{$\mid^{\text{d}}$\hss}\,
\lower4pt\hbox to 0pt{\hss$\smile$\hss}\ \ }}
\def\nindd{\mathop{\ \hbox to 
0pt{$\!\not{\mid}^{\text{\,d}}$\hss}\,
\lower4pt\hbox to 0pt{\hss$\smile$\hss}\ 
\ }}
\def\Mc{$\widetilde{\mathfrak M}_c$}

\def\M{\mathfrak M}

\def\p{\varphi}
\def\U{\mathcal{U}}
\def\Z{\tilde Z}
\def\C{\tilde C}
\def\Ut{U^\text{\th}}
\def\R{\mathbb R}
\def\proves{\vdash}
\def\tp{\mathrm{tp}}

\def\F{\mathbb F}

\newcommand{\ldim}{\operatorname{lin.dim}}
\def\Fix{\mathrm{Fix}}
\def\ale{\lesssim}

\def\aeq{\sim}

\def\acl{\mathrm{acl}}

\begin{document}

\title{Dimensional Groups and Fields}          
\author{Frank O. Wagner}
\date{\today}

\address{Universit\'e Lyon 1; CNRS; Institut Camille Jordan UMR 5208, 21
avenue
Claude Bernard, 69622 Villeurbanne-cedex, France}
\email{wagner@math.univ-lyon1.fr}

\begin{abstract} We shall define a general notion of dimension, and study groups and rings whose interpretable sets carry such a dimensio. In particular, we deduce chain conditions for groups, definability results for fields and domains, and show that a pseudofinite \Mc-group of finite positive dimension contains a finite-by-abelian subgroup of positive dimension, and a pseudofinite group of dimension 2 contains a soluble subgroup of dimension $2$.
\end{abstract}
\thanks{Partially supported by ValCoMo (ANR-13-BS01-0006)}

\maketitle

\section*{Introduction} In this paper, we shall define a general notion of dimension, and study groups and rings whose interpretable sets carry such a dimension. The aim is to unify results from stability and simplicity theory, $o$-minimality, and the study of pseudofinite structures (with dimensions induced by Lascar or SU-rank, $o$-minimal dimension, or the logarithm of the counting measure, see Example \ref{example}). In particular, we deduce chain conditions for groups, definability results for fields and domains, and show that pseudofinite groups contain big finite-by-abelian subgroups, and pseudofinite groups of dimension 2 contain big soluble subgroups.

\section{Dimension}
\defn\label{d:defdim} A theory $T$ is {\em dimensional} if there is a dimension function $\dim$ from the collection of all interpretable sets in models of $T$ to an ordered abelian group $\Gamma$ together with $\pm\infty$, 
satisfying for a formula $\p(x,y)$ and interpretable sets $X$ and $Y$:\begin{itemize}
\item\label{i:invariant} {\em Invariance:} If $a\equiv a'$ then $\dim(\p(x,a))=\dim(\p(x,a'))$.
\item\label{i:algebraic} {\em Algebraicity:} If $X$ is finite non-empty then $\dim(X)=0$, and $\dim(\emptyset)=-\infty$.
\item\label{i:union}{\em Union:} $\dim(X\cup Y)=\max\{\dim(X),\dim(Y)\}$.
\item\label{i:fibre} {\em Fibration:} If $f:X\to Y$ is a interpretable map such that $\dim(f^{-1}(y))\ge d$ for all $y\in Y$, then $\dim(X)\ge\dim(Y)+d$.
\end{itemize}
A dimension is {\em real} (or {\em archimedean}) if $\Gamma\le\mathbb R$, {\em discrete} if $\Gamma$ is discrete, and {\em integer} if it is discrete and real (thus $\Gamma\cong\mathbb Z$). Note that in a real or integer-dimensional theory sets of dimension $\infty$ are still allowed; a set (or a structure) $X$ is {\em finite-dimensional} if the dimension is integer and $\dim(X)<\infty$.\edefn
For a partial type $\pi$ put $\dim(\pi)=\inf\{\dim(\p):\pi\proves\p\}$, where the infimum is considered in a completion $\bar\Gamma$ of the ordered set $\Gamma\cup\{\infty\}$. Note that unless the dimension is real, there will be non-realized cuts in $\Gamma$ (for instance the $\sup$ of a proper convex subgroup), and the extension of the group-operation by continuity is only well-defined on the semigroup $\bar\Gamma^+=\{\gamma\in\bar\Gamma:\gamma\ge0\}$. We 
write $\dim(a/B)$ for $\dim(\tp(a/B))$.  Then for types dimension is invariant, algebraic and satisfies union, but need not satisfy fibration. By the union condition any partial type can be completed to a type of the same dimension.

\bem\label{r:defdim} There are some additional axioms and variants one might also consider:\begin{itemize}
\item\label{i:fine}{\em Finesse:} If $\dim(X)=0$ then $X$ is finite.
\item\label{i:product} {\em Product:} $\dim(X\times Y)=\dim(X)+\dim(Y)$.
\item\label{i:sfibre}{\em Strong fibration:} If $f:X\to Y$ is a interpretable map such that $\dim(f^{-1}(y))= d$ for all $y\in Y$, then $\dim(X)=\dim(Y)+d$.
\item\label{i:wfibre}{\em Weak fibration:} If $f:X\to Y$ is a interpretable map such that $\dim(f^{-1}(y))=d$ for all $y\in Y$, then $\dim(X)\ge\dim(Y)+d$.
\item{\em Lower fibration:} If $f:X\to Y$ is a interpretable map such that $\dim(f^{-1}(y))\le d$ for all $y\in Y$, then $\dim(X)\le\dim(Y)+d$.
\item\label{i:definability}{\em Definability:} If $\dim(\p(x,a))=d$ then there is a formula $\theta_{\p,d}\in\tp(a)$ such that $\dim(\p(x,a'))=d$ for all $a'\models\theta_{\p,d}$.
\item\label{i:+definability}{\em Semidefinability:} If $\dim(\p(x,a))>d$ then there is a formula $\theta_{\p,d}\in\tp(a)$ such that $\dim(\p(x,a'))> d$ for all $a'\models\theta_{\p,d}$.
\item\label{i:additive}{\em Additivity:} $\dim(a,b/A)=\dim(a/A,b)+\dim(b/A)$.
\item\label{i:+additive}{\em Semiadditivity:} $\dim(a,b/A)\ge\dim(a/A,b)+\dim(b/A)$.
\end{itemize}
We call a dimension {\em coarse} if it is not fine, i.e.\ does not satisfy finesse.
\ebem

Clearly both fibration and strong fibration imply weak fibration, strong fibration implies product, fibration and lower fibration imply strong fibration, any kind of definability implies invariance, definability plus strong fibration imply additivity, and semidefinability plus fibration imply semiadditivity.
\bem\begin{itemize}\item Invariance is equivalent to type-definability, i.e.\ definability where $\theta_{\p,d}$ is a partial type (for instance $\tp(a)$ itself will do).
\item We have chosen weak inequalities for (weak) fibration, as this behaves better under limits. On the other hand, we took strict inequalities for semidefinability, as this seems easier to achieve in examples.
\item One might also define lower semidefinability where the inequalities are reversed. This seems less useful, though.
\item Note that fibration yields the inequality $\dim(X\times Y)\ge\dim(X)+\dim(Y)$.
\item If definability holds, then by compactness $\dim$ only takes finitely many values in any given sort or arity.
\end{itemize}
\ebem
\bem\label{lfibration} Additivity implies fibration and lower fibration, whence strong fibration.\ebem
\bew Let $f:X\to Y$ be a definable map, and note that $\dim(f(x)/x)=0$ for all $x\in X$, whence $$\dim(x)=\dim(x)+\dim(f(x)/x)=\dim(x,f(x))=\dim(x/f(x))+\dim(f(x)).$$
Suppose first that $\dim(f^{-1}(y))\ge d$ for all $y\in Y$. Take $b\in Y$ with $\dim(b)=\dim(Y)$, and choose $a\in f^{-1}(b)$ with $\dim(a/b)\ge d$.  Then
$$\dim(X)\ge\dim(a)=\dim(a/b)+\dim(b)\ge d+\dim(Y).$$
Suppose now that $\dim(f^{-1}(y))\le d$ for all $y\in Y$. Take $a\in X$ with $\dim(a)=\dim(X)$, and put $b=f(a)$. Then
$$\dim(X)=\dim(a)=\dim(a/b)+\dim(b)\le d+\dim(Y).\qedhere$$
\ebew

\bsp\label{example} Examples for integer dimensions with lower/strong fibration include:\begin{enumerate}
\item\label{e:item} Finite Lascar rank, $SU$-rank or $U^{\mbox{\tiny\th}}$-rank on {\em formulas}, possibly localised at some 
$\emptyset$-invariant family of types;
\item For any ordinal $\alpha$, the coefficient of $\omega^\alpha$ in one of 
the ordinal-valued ranks in (\ref{e:item}) above (when written on Cantor normal form).
\end{enumerate}
In these examples, for an interpretable set $X$ we let  $\dim(X)$ be the maximum of the finite ranks (Lascar, $SU$ ou $U^{\mbox{\tiny\th}}$) of the types extending $x\in X$, and $\dim(X)=\infty$ if no such maximum exists. In particular, in the case of finite Lascar rank, the dimension of a formula is just equal to its Shelah rank. Note that in general the dimension of a type need not be equal to the rank of the type even when the dimension is finite, as witnessed by the standard example where Lascar and Shelah rank are finite and different.\footnote{A structure with disjoint unary predicates $\{P_i:i<\omega\}$, and on every $P_i$ an equivalence relation $E_i$ with infinitely many infinite classes. The type $\{\neg P_i:i<\omega\}$ has Lascar rank $1$ but Shelah rank~$2$.} Nevertheless, the Lascar inequalities for types are sufficient to show that fibration and lower fibration hold (whence strong fibration as well). Note that in example (1), not localised, the dimension is fine; if the rank is finite, so is the dimension.
\begin{enumerate}\item[(4)] $o$-minimal dimension on a densely ordered set. This dimension is finite, and satisfies all additionnal properties of Remark \ref{r:defdim}.
\end{enumerate}
An example for a real additive dimension is\begin{enumerate}
\item[(5)] coarse pseudofinite dimension (in some expansion by cardinality comparison quantifiers), defined for a definable set $X=\prod_i X_i/\U$ in a pseudofinite structure $\prod_i\M_i/\U$ as $\log(\prod|X_i|/\U)+C\in(\prod_i\R/\U)/C$, where $C$ is a convex additive subgroup containing the integers in the non-standard real field $\prod_i\R/\U$, and $|X_i|$ is the cardinality of $X_i$, see \cite{Hru12, Hru13, HW08}. If $C$ is the convex hull of the integers, this is {\em fine} pseudofinite dimension.
\end{enumerate}
Here additivity holds by \cite[Lemma 2.10]{Hru13}, and in fact all additionnal properties of Remark \ref{r:defdim} hold except for definability and finesse (for semidefinability one needs cardinality comparison quantifiers). Clearly fine pseudofinite dimension is fine in the sense of Remark~\ref{r:defdim}.
\ebsp

\bem If $\dim$ is a $\Gamma$-valued dimension and $\gamma\in\Gamma$, put 
$$\begin{aligned}\Gamma_\gamma&=\{\gamma'\in\Gamma:|\gamma'|\le n|\gamma|\mbox{ for some $n<\omega$}\},\quad\mbox{and}\\
\Gamma_{\gamma^-}&=\{\gamma'\in\Gamma:n|\gamma'|< |\gamma|\mbox{ for all $n<\omega$}\}.
\end{aligned}$$
Then $\Gamma_{\gamma^-}<\Gamma_\gamma\le\Gamma$ are subgroups, and there is a unique additive monomorphism $\sigma:\Gamma_\gamma/\Gamma_{\gamma^-}\to\mathbb R$ with $\sigma(\gamma)=1$. The function $\dim_\gamma$ defined by
$$\dim_\gamma(X)=\left\{\begin{array}{ll}\sigma(\dim(X))\quad&\mbox{if $\dim(X)\in\Gamma_\gamma$}\\
\infty&\mbox{otherwise}\end{array}\right.$$
is a real dimension, the {\em localisation} of $\dim$ at $\gamma$. If $\gamma=\dim(Y)$ for some definable set $Y$, we also write $\dim_Y$ instead of $\dim_\gamma$.\ebem

\lem\label{l:weak} Just assuming weak fibration, dimension is invariant under definable bijections. Under strong fibration, dimension is invariant under finite-to-finite correspondences.\elem
\bew Let $f:X\to Y$ be a definable bijection. Then $\dim(f^{-1}(y))=0$ for all $y\in Y$, so $\dim(X)\ge\dim(Y)$. Considering the definable bijection $f^{-1}:Y\to X$, we obtain $\dim(Y)\ge\dim(X)$, whence equality.

Now assume strong fibration, and let $R\subset X\times Y$ be a finite-to-finite correspondence with non-empty fibres of size at most $n$. For $x\in X$ put $Y_x=\{y\in Y:(x,y)\in R\}$, and consider $Z=\{Y_x:x\in X\}$. Then the map $x\mapsto Y_x$ from $X$ to $Z$ has finite fibres bounded by $n$, so $\dim(X)=\dim(Z)$ by strong fibration. Next, consider $Z'=\{(y_1,\ldots,y_n)\in Y^n:\{y_1,\ldots,y_n\}\in Z\}$. Then the map $(y_1,\ldots,y_n)\mapsto \{y_1,\ldots,y_n\}$ from $Z'$ to $Z$ has non-empty fibres of size at most $n!$, so $\dim(Z')=\dim(Z)$ by strong fibration. Finally, the projection to the first coordinate from $Z'$ to $Y$ has non-empty fibres of size at most $n^n$: for any $y_1\in Y$ there are at most $n$ different $x\in X$ with $(x,y)\in R$, and for each $x$ there are at most $n^{n-1}$ choices for the coordinates $y_2,\ldots,y_n$. Thus $\dim(Z')=\dim(Y)$.\ebew

\defn Let $X$ and $Y$ be type-definable sets. We say that $X$ is {\em broad} if $0<\dim(X)<\infty$. If $X$ is broad, we say that $Y$ is {\em $X$-broad} if $\frac1n\dim(X)\le\dim(Y)\le n\dim(X)$ for some $n<\omega$.

A definable set $X$ is {\em negligible} if $\dim(X)=0$; a type-definable set is {\em negligible} if it is contained in some negligible definable set.

A type-definable subset $Y\subseteq X$ is {\em wide} in $X$ if for every definable superset $\bar Y\supseteq Y$ there is a definable superset $\bar X\supseteq X$ with $\dim(\bar Y)\ge\dim(\bar X)$; in particular a definable subset $Y$ of a definable set $X$ is wide iff $\dim(Y)=\dim(X)$. An element $x\in X$ is {\em wide}/{\em broad} over some parameters $A$ if $\tp(x/A)$ is.\edefn

\bem\begin{itemize}\item $X$ is $Y$-broad iff $Y$ is $X$-broad.
\item If $X$ is broad, $Y$ is $X$-broad iff $Y$ is broad for the localized dimension $\dim_X$.
\item If the dimension is real and $X$ is broad, then $Y$ is broad iff it is $X$-broad.
\item In fine dimension, a negligible type-definable set is finite.
\item In fine, discrete dimension, a type-definable set of dimension $0$ is finite.
\end{itemize}\ebem

\section{Dimensional Groups}
\lem\label{generics} Let $G$ be a type-definable group over $A$ in a dimensional theory, and $g,h\in G$. If $h$ is wide in $G$ over $A,g$, then $gh$ and $hg$ are wide in $G$ over $A,g$.\elem
\bew Let $X$ be an $A,g$-definable set containing $gh$. Since $G$ is a type-definable group, restricting $X$ we may assume that $x\mapsto g^{-1}x$ is a bijection between $X$ and $g^{-1}X$. As dimension is invariant under definable bijections, $\dim(X)=\dim(g^{-1}X)$. But $h\in g^{-1}X$, so $gh$ is wide over $A,g$. The proof for $hg$ is similar.\ebew

In a dimensional theory, we need not have fibration for definable maps between type-definable sets. The situation is different for group homomorphisms.
\lem\label{l:dhomo}
In a dimensional theory, let $G$ and $H$ be type-definable groups and $f: G \to H$ a definable surjective homomorphism. Then $\dim(G) \ge\dim(\ker f) + \dim(H)$.
\elem
\bew If $\dim(G)=\infty$ this is clear. Otherwise, consider a definable $X_0\supseteq G$. Reducing $X_0$, if necessary, we may assume by compactness that there is a definable map $\bar{f}$ extending $f$ with domain $X_0$. Again by compactness there is a definable $X_1=X_1^{-1}\supseteq G$ with $X_1^2\subseteq X_0$, and such that $\bar{f}(xx')=\bar{f}(x)\bar{f}(x')$ for all $x,x'\in X_1$. Put $Y=\bar f(X_1)$ and $X=\bar f^{-1}(Y)\subseteq \mathop{\mathrm{dom}}(\bar f)=X_0$. Then for $y\in Y$ we have
$$\bar{f}^{-1}(y)\supseteq (\bar{f}^{-1}(y)\cap X_1)\ker f,$$
so $\dim(\bar f^{-1}(y))\ge\dim(\ker f)$. Since $H\subseteq Y$ it follows that
\[\dim(X_0)\ge\dim(X)\ge\dim(\ker f)+\dim(Y)\ge\dim(\ker f)+\dim(H).\]
Thus $\dim(G)=\inf\{\dim(X_0):X_0\supseteq G\mbox{ definable}\}\ge\dim(\ker f)+\dim(H)$.\ebew

We now turn to chain conditions.
\satzli\label{l:dimensiondrops-d}
In a fine dimensional theory, let $H$ be a definable subgroup of infinite index in a type-definable group $G$ with $\dim(G)<\infty$. Then $\dim(H)<\dim(G)$.
\esatzli
\bew Since $G$ and $H$ are definable, the quotient space $G/H$ is interpretable. Thus $\dim(G/H)$ is well-defined, and strictly positive by finesse. The map $G\to G/H$ has fibres of dimension $\dim(H)$. Hence by  fibration,
\[\dim(G)\ge\dim(H)+\dim(G/H)>\dim(H).\qedhere\]
\ebew
This immediately yields:
\kor\label{c:chain-d}
In a fine integer-dimensional theory there is no infinite descending chain $(G_i:i<\omega)$ of definable groups of finite dimension, each of infinite index in its predecessor.\qed\ekor

Since we have not required the dimension to be defined on quotients by \emph{type-definable} equivalence relations, Proposition~\ref{l:dimensiondrops-d} and Corollary~\ref{c:chain-d} may fail for type-definable subgroups. This happens for instance in $o$-minimal theories: the additive subgroup of infinitesimals has the same dimension as the ambient non-standard real field, but infinite, and even unbounded, index. In fact, we can have arbitrarily long (infinite) chains of more and more infinitesimal type-definable subgroups. However, we shall show next that this does not happen for a particular kind of type-definable groups.

\defn
We call a group \emph{(relatively) $\bigwedge$-definable} if it is an intersection of (relatively) definable groups. A ring is (relatively) $\bigwedge$-definable if it is so as an additive group.
\edefn
As we have not defined dimension on arbitrary hyperimaginaries (quotients modulo type-definable equivalence relations), we do not have a dimension on the quotient of two type-definable groups. However, for a quotient of a type-definable group $G=\bigwedge_{i\in I} X_i$ by a relatively $\bigwedge$-definable subgroup $H=\bigwedge_{j\in J}H_j$ (where the $X_i$ are definable sets, the $H_j$ are definable subgroups, and we take both to be closed under finite intersections for ease of notation), we can put 
$$\dim(G/H)=\sup_{j\in J}\dim(G/H_j)=\sup_{j\in J}\inf_{i\in I}\dim(X_i/H_j)\in\bar\Gamma.$$

\bem\begin{itemize}
\item
In an $\omega$-stable theory (and in particular in a theory of finite Morley rank), definable, $\bigwedge$-definable and type-definable groups coincide.
\item
In a stable or supersimple theory, type-definable groups are $\bigwedge$-definable.
\item
In an $o$-minimal theory (as in any theory with descending chain condition on definable groups), $\bigwedge$-definable subgroups are definable.
\item
Even in the fine integer-dimensional context, the three classes need not coincide. For instance:
\begin{itemize}
\item
The connected component of $\mathbb{Z}$ (in a saturated model) is a $\bigwedge$-definable group which is not definable, in a superstable theory of Lascar rank $1$.
\item
In a non-standard real field the infinitesimals form a type-definable additive group which is not $\bigwedge$-definable, in an $o$-minimal theory.
\end{itemize}
\item If $(G_i:i<\omega)$ is a chain of definable groups with $\dim(G_i/G_{i+1})=1$ for all $i<\omega$, our definition yields $\dim(G/G)=0$ as expected. Had we exchanged the limits, we would get $\inf_i\sup_j\dim(G_i/G_j)=\infty$, which is clearly wrong.
\end{itemize}
\ebem

We first check that $\bigwedge$-definability behaves well with respect to relative definability.

\lem\label{l:relativeiswedgedefinable}
Let $G$ be a $\bigwedge$-definable group, and $H$ a relatively definable subgroup. Then $H$ is $\bigwedge$-definable. In fact, there is a definable group $H^*$ with $H=G\cap H^*$.
\elem
\bew
Suppose $G = \bigwedge_{i\in I} G_i$, where the $G_i$ are definable groups, and $H = G\cap X$, where $X$ is a definable set.
Then $x,y\in G\cap X$ implies $x^{-1}y \in X$. By compactness there is some finite $I_0 \subseteq I$ such that $x,y \in \bigcap_{i\in I_0} G_i \cap X$ implies $x^{-1}y \in X$. As the $G_i$ are groups, we also have $x^{-1}y \in \bigcap_{i \in I_0} G_i$, and $Y = \bigcap_{i\in I_0} G_i \cap X$ is a definable group with $H = G \cap Y$.
\ebew

\satzli\label{l:dimensiondrops}
In a dimensional theory, let $H$ be a $\bigwedge$-definable subgroup in a type-definable group $G$. Then $\dim(G)\ge\dim(H)+\dim(G/H)$. In particular, if the dimension is fine and discrete, or if the dimension is fine and $G$ is definable, then $\dim(H)<\dim(G)$ iff $H$ has unbounded index in $G$.
\esatzli
\bew Suppose $H=\bigwedge_i H_i$, where each $H_i$ is a relatively definable subgroup of $G$ and the system $(H_i)_i$ is closed under finite intersection. Then for every $i$ the projection $G\to G/H_i$ has fibres of dimension $\dim H_i$, whence by Lemma \ref{l:dhomo}
$$\dim(G)\ge\dim(H_i)+\dim(G/H_i)\ge\dim(H)+\dim(G/H_i).$$
Therefore
$$\begin{aligned}\dim(G)&\ge\sup_i[\dim(H)+\dim(G/H_i)]\\
&=\dim(H)+\sup_i\dim(G/H_i)=\dim(H)+\dim(G/H).\end{aligned}$$
In fine dimension, if the dimension is discrete or $G$ is definable, then 
$$\begin{aligned}G/H\mbox{ is unbounded}&\Leftrightarrow G/H_i\mbox{ infinite for some $i$}\\
&\Leftrightarrow \dim(G/H_i)>0\mbox{ for some $i$}\Leftrightarrow\dim(G/H)>0.\end{aligned}$$
Thus under either condition, $\dim(G)>\dim(H)$ iff $G/H$ is unbounded.\ebew

\kor\label{c:chain}
In a fine integer-dimensional theory there is no infinite descending chain of relatively $\bigwedge$-definable groups of finite dimension, each of unbounded index in its predecessor.
\ekor
\bew Let $(G_i:i<\omega)$ be such a chain. As the dimension is finite and discrete, $\dim(G_i)>\dim(G_{i+1})$ by Proposition \ref{l:dimensiondrops}. But there is no infinite descending sequence of positive integers.
\ebew

\bem In a coarse integer-dimensional theory we can still conclude that there is no infinite descending  chain of relatively $\bigwedge$-definable groups of finite dimension, with non-negligible successive quotients.\ebem

\section{Fields and Domains}
\subsection{Skew fields} Let us note first that fields have better definability properties than groups. Our first result generalizes the well-known fact that a supersimple type-definable field is in fact definable.

\satzli\label{l:fielddef}
In a real-dimensional theory, a type-definable broad skew field $K$ is definable.
\esatzli
\bew
Suppose $K=\bigcap_{i\in I} X_i$, where $(X_i: i\in I)$ is a system of definable sets closed under finite intersections. As $0<\dim(K)<\infty$ we may also assume $\dim(X_i)<2\dim(K)$ for all $i\in I$. By compactness, we may further suppose that there is a minimal element $0\in I$ such that addition and multiplication are defined, commutative, associative and distributive on $X_0$ (but may take values outside), and for every $i>0$ all non-zero elements in $X_i$ have an additive and a multiplicative inverse in $X_i$.

By compactness there is $i\in I$ such that $X_i\cdot (X_i-X_i)+(X_i-X_i)\subseteq X_0$. But then for any $X_j\subseteq X_i$ and $g\in X_i$ we have $g X_j + X_j\subseteq X_0$. Moreover, if $g(X_j-X_j)\cap(X_j-X_j)=\{0\}$, then the map from $X_j^2\to X_0$ given by $(x,y)\mapsto gx+y$ is injective, contradicting 
$$\dim(X_j^2) \ge 2 \dim(X_j) \ge 2 \dim(K) > \dim(X_0).$$
Thus $g\in (X_j-X_j)\cdot\big((X_j-X_j)\setminus\{0\}\big)^{-1}$, and
\[X_i\subseteq\bigcap_{X_j\subseteq X_i} (X_j-X_j)\cdot\big((X_j-X_j)\setminus\{0\}\big)^{-1}.\]
But by compactness for every $k\in I$ there is some $j\in I$ such that
\[(X_j-X_j)\cdot\big((X_j-X_j)\setminus\{0\}\big) ^{-1}\subseteq X_k.\]
Thus
\[K\subseteq X_i
\subseteq \bigcap_j (X_j-X_j)\cdot\big((X_j-X_j)\setminus\{0\}\big)^{-1}
\subseteq  \bigcap_k X_k \subseteq K,\]
and $K=X_i$ is definable.
\ebew
\bem Note that the hypothesis of Proposition \ref{l:fielddef} requires that $\dim(K)$ is bounded away from $0$.\ebem

\satzli\label{l:Vdef}
In a real-dimensional theory, let $K$ be a non-negligeable definable skew field and $V$ a non-trivial type-definable $K$-vector space with $\dim(V)<\infty$. Then $\dim(V)\ge\ldim_K(V)\,\dim(K)$. In particular $K$ and $V$ are broad, $\ldim_K(V)$ is finite and $V$ is definable.
\esatzli
\bew Clearly $\dim(K^n)\ge n\,\dim(K)$. As dimension is preserved under definable bijection, $\dim(K)>0$ and $\dim(V)<\infty$, we get $$\infty>\dim(V)\ge\ldim_KV\, \dim(K)\ge\dim(K)>0.$$
So $K$ is broad, $\ldim_KV$ is finite, and $V$ is definable as $\sum_i Ke_i$, for a $K$-basis $(e_i)_i$ of $V$.
\ebew
\bem If the dimension satisfies {\em product}, then $\dim V=\ldim_KV\,\dim K$.\ebem

\satzli\label{l:centre}
Let $K$ be a definable broad skew field in a fine real-dimensional theory. Then $K$ has finite dimension over its centre.\esatzli
\bew
If $K^\times$ has finite exponent, then $K$ satisfies a polynomial identity and has finite dimension over its centre by Kaplansky's PI-Theorem \cite{Ka48}. In fact, the centre is a field of finite exponent, whence finite, and so is $K$.

Otherwise, there is an element $a\in K$ of infinite order, and $Z(C_K(a))$ is an infinite definable commutative subfield. But then $K$ has finite dimension over the infinite definable subfield $Z(C_K(a))$, whence finite dimension over its centre.
\ebew

\subsection{Domains} We now move to domains, generalizing both \cite[Theorem~2.2]{Krup10} and \cite[Exercises 3.5 and 3.6]{Hru13}. Recall that localising a non-commutative domain is not always possible: not only may right and left fractions differ, it is not a priori possible to multiply say two right fractions in a consistent way.

The right/left \emph{Ore condition} in a domain $R$ asks that $rR\cap r'R\not=\{0\}$ for all non-zero $r,r'\in R$ (or $Rr\cap Rr'\not=\{0\}$, respectively). If it holds then there is a right (resp., left) fraction skew field; if both hold they give rise to the same skew field. A domain is a left/right {\em Ore domain} if it satisfies the left/right Ore condition.

\satzli\label{l:ore}
In a real-dimensional theory, let $R$ be an invariant domain, and suppose that there is a constant $d<\infty$ such that $\dim(X)\le d$ for all type-definable $X\subseteq R$, and $\dim(X)>0$ for some such $X$. Then $R$ is right and left Ore, and its skew field $K$ of (right or left) fractions is definable.
\esatzli
\bew As $R$ is invariant, it is a union of type-definable sets (over some set of parameters).
The assumptions imply that there is some type-definable $X\subseteq R$ with $2\dim(X)>\dim(Y)$ for any type-definable $Y\subseteq R$. For any two non-zero $r,r'\in R$ consider the map
\[\begin{array}{ccc}
X^2 & \to & rX + r'X\\
(x,y)& \mapsto & rx+r'y.
\end{array}\]
Suppose it is injective. As $rX + r'X$ is a type-definable subset of $R$, the choice of $X$ implies $\dim(rX+r'X)<2\dim(X)$. Hence there is a definable superset $Y$ of $X$ with $\dim(rY+r'Y)<2\dim(X)$, and such that the map $(x,y)\mapsto rx+r'y$ is bijective between $Y^2$ and $rY+r'Y$. Then
$$2\dim(X)>\dim(rY+r'Y)=\dim(Y^2)\ge 2\dim(Y)\ge 2\dim(X)$$
(where the equality holds by Lemma \ref{l:weak} and the first weak inequality by Remark \ref{r:defdim}), a contradiction.

Thus there are $(y,y') \neq (z,z')$ with $ry+r'y' = rz+r'z'$.
Then $r(y-z)=r'(z'-y') \neq 0$, so $R$ is right Ore. Similarly, $R$ is left Ore. It follows that its field of left fractions is equal to its field of right fractions, and equal to $(X-X)/(X-X)^\times$, which is type-definable. Hence $K$ is definable by Proposition~\ref{l:fielddef}.
\ebew

\satzli\label{l:dimdomain}
In a fine dimensional theory satisfying product, a broad $\bigwedge$-definable (non-commu\-tative, non-unitary) left or right Ore domain $R$ with definable fraction field $K$ is already a skew field.
\esatzli
\bew We may assume that $R$ is infinite by Wedderburn's Theorem. Consider a definable additive supergroup $A$ of $R$ such that inversion and multiplication is well-defined and associative on $A^{\pm1}$, and without zero-divisors. Choose a definable $B$ with $R\le B\le A$ (as additive groups) and with $B+B^2\subseteq A$.

Consider the map $f:(x,y)\mapsto xy^{-1}$ on $A\times A^\times$, put $Y=B{B^\times}^{-1}\supseteq K$ and $X=f^{-1}(Y)\subseteq A\times A^\times$. For every $y\in Y$ we have
$$\dim(f^{-1}(y))\ge\dim((b,b')B^\times)=\dim(B)$$
for any $(b,b')\in B\times B^\times$ with $y=bb'^{-1}$.

Suppose $\dim(R)<\dim(K)$. Then $2\dim(R)<\dim(K)+\dim(R)$, so by $\bigwedge$-definability of $R$, we can choose $A$ (and $B$) such that $2\dim(A)<\dim(K)+\dim(R)$. But product and fibration for $f$ yield
$$2\dim(A)=\dim(A^2)\ge\dim(X)\ge\dim(B)+\dim(Y)\ge\dim(R)+\dim(K)>2\dim(A),$$
a contradiction. It follows that $\dim(R)=\dim(B)=\dim(K)$. 

As the dimension is fine and $K$ definable, $B$ has finite index in $(K,+)$ by Proposition \ref{l:dimensiondrops-d}. Hence for any $r\in R^\times$ there are natural numbers $m,n$ with $n+1<m$ such that $r^{-m}$ and $r^{-n}$ lie in the same coset modulo $B$. Thus there is $b\in B$ with $r^{-m}=r^{-n}+b$, and 
$$r^{-1}=r^{m-1}r^{-m}=r^{m-1}(r^{-n}+b)=r^{m-n-1}+r^{m-1}b\in R+RB\subseteq B+B^2\subseteq A.$$
As this also holds for any definable subgroup of $A$ containing $R$, we get $r^{-1}\in R$ by $\bigwedge$-definability, whence $R=K$.\ebew

\kor\label{c:skew} In a fine real-dimensional theory, a $\bigwedge$-definable broad (non-commu\-tative, non-unitary) domain is a definable skew field.\qed\ekor

\bem The infinitesimals in a non-standard real closed field show that Proposition \ref{l:dimdomain} and Corollary  \ref{c:skew} may fail for type-definable domains in an $o$-minimal (and hence fine integer-dimensional) theory.\ebem

\kor\label{c:prelinearisation}
In a real-dimensional theory, let $A$ be a type-definable abelian group with $\dim(A)<\infty$. Suppose that there is an invariant set $X$ generating a domain $R$ of definable automorphisms of $A$, and such that $\dim(Y)>$ for some type-definable $Y\subseteq X$. Then the skew field of fractions $K$ of $R$ exists and is definable; $A$ is a definable $K$-vector space of finite linear dimension.
\ekor
\bew
Fix any non-zero $a \in A$. Since $R$ acts by automorphisms, the evaluation map $R \to R\cdot a$ is injective, which bounds $\dim (Z) = \dim(Z \cdot a)$ by $\dim(A)$ for any type-definable subset $Z\subseteq R$. Now apply  Proposition~\ref{l:ore}, and note that the fraction field $K$ acts on $A$ since $R^\times$ acts by automorphisms. We finish by Proposition \ref{l:Vdef}.
\ebew

Notice that if $R$ is commutative, so is $K$.

\bem In Corollary~\ref{c:prelinearisation} the hypothesis that $\dim(Y)>0$ for some type-definable $Y\subseteq X$ is necessary: Just consider an infinite field $K$ in the language of modules consisting of addition and unary functions $\lambda_r$ for scalar multiplication by $r\in K$, for every $r\in K$. This is an abelian structure which does not interpret an infinite field.\ebem

\subsection{Automorphisms}
As opposed to the stable case, an infinite field with integer fine dimension need not be $\emptyset$-connected and may have $\emptyset$-definable additive or multiplicative subgroups of finite index: any pseudo-finite field will serve as an example.

\lem
In a real-dimensional theory, a definable endomorphism $\varphi$ of a definable broad skew field $K$ is either $0$ or a genuine skew field automorphism.
\elem
\bew If $\varphi$ is not zero, it is injective, and $0<\dim(\varphi(K)) = \dim(K)<\infty$. Hence the degee $[K:\varphi(K)]=1$ and $K=\varphi(K)$.\ebew

It follows that a dimensional broad (commutative) field is perfect.

Next, we show that a $\emptyset$-connected definable broad field in a fine real-dimensional theory does not admit an infinite type-definable family of automorphisms. This generalises \cite[Theorem~8.3]{BNGroups}. Recall that an element of a field is {\em absolutely algebraic} if it is the root of a non-zero polynomial with coefficients in $\mathbb Z$; the set of absolutely algebraic elements forms a subfield containing the prime field.

\lem\label{l:absalg}In a fine dimensional theory, a $\emptyset$-definable broad additively or multiplicatively $\emptyset$-connected field $K$ contains infinitely many absolutely algebraic elements.\elem
\bew In characteristic zero this is clear. So suppose the field $K$ has characteristic $p>0$ and is additively $\emptyset$-connected; consider the $\emptyset$-definable additive endomorphism $\phi:x \mapsto x^p - x$ with finite kernel $\F_p$. If there are only finitely many absolutely algebraic elements, then there is some $n<\omega$ such that $\phi^n(K)\cap\F_p=\{0\}$. So $\phi:\phi^n(K)\to\phi^{n+1}(K)$ is injective, whence $\dim(\phi^{n+1}(K))=\dim(\phi^n(K))$; as the dimension is fine, $\phi^{n+1}(K)$ has finite index in $\phi^n(K)$. But $\emptyset$-connectivity of $K$ implies that of $\phi^n(K)$, and $\phi^n(K)=\phi^{n+1}(K)=\phi^{2n}(K)$. But then $K=\phi^n(K)\oplus\ker\phi^n$, so $\phi^n(K)$ has finite index in $K$ as an additive subgroup. Again by connectivity $\phi$ is surjective. Hence $K$ has no Artin-Schreier extension. But the subfield of absolutely algebraic elements is relatively algebraically closed in $K$; if it were finite, it would have an extension of degree $p$, which would yield an Artin-Schreier extension of $K$.

If $K$ is multiplicatively $\emptyset$-connected, we use the multiplicative endomomorphism $\phi:x\mapsto x^q$ for big prime $q$ such that $K$ has no primitive $q$-th root of unity, and argue with Kummer extensions.\ebew

\satzli\label{l:fieldautomorphisms}
In a fine real-dimensional theory, let $K$ be a definable broad field with infinitely many absolutely algebraic elements. Then there is no infinite type-definable family of definable automorphisms of $K$.
\esatzli
\bew
Suppose $\Phi$ is such a family. As $K$ is definable, we may assume that $\Phi$ is also definable (being an automorphism of $K$ is a definable property).
If char$(K)=0$, put $K_0=\bigcap_{\sigma\in\Phi}\Fix(\sigma)$. Then $K_0$ is definable and infinite, whence broad, and $[K:K_0]$ is finite by Proposition \ref{l:Vdef}. Thus $\Phi$ is a subfamily of Gal$(K/K_0)$, which is finite.

Now suppose char$(K)=p>0$; we put 
$$\Psi=\{\tau^{-1}\sigma:\sigma,\tau\in\Phi\}.$$
If $a\in K$ is absolutely algebraic, it has only finitely many images under $\Phi$, and there are infinitely many automorphisms in $\Psi$ fixing $a$. By compacteness there is some $a$ of infinite multiplicative order fixed by infinitely many automorphisms of $\Psi$. Put $\Psi_a=\{\sigma\in\Psi:\sigma(a)=a\}$ and
$$K_0=\bigcap_{\sigma\in\Psi_a}\Fix(\sigma).$$
Then $K_0$ is definable and infinite, whence broad, and $\Psi_a$ injects into Gal$(K/K_0)$ which is finite.

In both cases we obtain a contradiction, so no infinite type-definable family of automorphisms of $K$ can exist.\ebew

\kor\label{c:fieldautomorphisms}In a fine real-dimensional theory, a definable automorphism of a $\emptyset$-definable, broad, additively or multiplicatively $\emptyset$-connected field is $\acl^{eq}(\emptyset)$-definable.
\ekor
\bew The field has infinitely many absolutely algebraic elements by Lemma \ref{l:absalg}. If there were an $a$-definable automorphism $\sigma_a$ not definable over $\acl^{eq}(\emptyset)$, there would be an infinite type-definable family 
$\{\sigma_{a'}:a'\models\tp(a)\}$ of definable automorphisms, contradicting  Proposition \ref{l:fieldautomorphisms}.
\ebew

\section{Pseudofinite dimensional groups}
\satzli\label{centralizer} Let $G$ be a broad pseudofinite dimensional group with strong fibration. Then there is an element $g\in G\setminus\{1\}$ with $G$-broad centralizer. More precisely, there is $g\in G\setminus\{1\}$ with $\dim(C_G(g))\ge\frac13\dim(G)$.\esatzli
\bew Suppose first that $G$ has no involution.
If $G\equiv\prod_I G_i/\U$ for some family $(G_i)_I$ of finite groups and some 
non-principal ultrafilter $\U$, then $G_i$ has no involution for almost all 
$i\in I$, and is soluble by the Feit-Thompson theorem. So for almost all $i\in I$ there is $g_i\in G_i\setminus\{1\}$ such that 
$\langle g_i^{G_i}\rangle$ is commutative. Put $g=[g_i]_I\in 
G\setminus\{1\}$. Then $\langle 
g^G\rangle$ is commutative and $g^G\subseteq C_G(g)$. As $g^G$ is in definable 
bijection with $G/C_G(g)$, we have 
$$\dim(C_G(g))\ge\dim(g^G)=\dim(G/C_G(g))=\dim(G)-\dim(C_G(g)).$$
In particular $\dim(C_G(g))\ge\frac12\dim(G)$.

Now suppose $G$ has an involution $i$, but all centralizers of non-trivial elements have dimension $<\frac13\dim(G)$. Then 
$$\dim(i^G)=\dim(G/C_G(i))=\dim(G)-\dim(C_G(i))>\frac23\dim(G).$$
For $h\in G\setminus\{1\}$ put $H_h=\{x\in G:h^x=h^{\pm1}\}$. Then $H_h$ is an $h$-definable subgroup of $G$, 
and $C_G(h)$ has index two in $H_h$, so $\dim(H_h)=\dim(C_G(h))<\frac13\dim(G)$. Moreover, if $j\in i^G$ and $h\in j^Gj$, then $j\in H_h$. Now by strong fibration
$$\dim(\{(j,h)\in G\times G: j\in H_h\})=\dim(G)+\dim(H_h)<\frac43\dim(G).$$
On the other hand by fibration
$$\begin{aligned}\dim(\{(j,h)\in G\times G: j\in H_h\})&\ge \dim(\{(j,h)\in i^G\times G: h\in j^Gj\})\\
&\ge2\dim(i^G)>\frac43\dim(G),\end{aligned}$$
as $\dim(j^Gj)=\dim(j^G)=\dim(i^G)>\frac23\dim(G)$ for all $j\in i^G$. This contradiction finishes the proof.
\ebew

We recall from \cite{He} the definition of \Mc, the centraliser condition up to finite index.
\defn A group $G$ satisfies the \Mc-condition if there is $n<\omega$ such that there are no $(g_i:i<n)$ in $G$ such that $|C_G(g_j:j<i):C_G(g_j:j\le i)|\ge n$ for all $i<n$. In other words, in a saturated model there is no infinite chain of centralisers $C_G(g_j:j<i)$ for $i<\omega$, each of infinite index in its predecessor.
\edefn
Examples for \Mc-groups include all groups definable in a simple theory. Note that a subgroup of an \Mc-group is again \Mc.
\lem\label{byZ} Let $G$ be an \Mc-group and $Z$ a finite central subgroup. Then $G/Z$ is an \Mc-group.\elem
\bew As $Z$ is central, $x\mapsto [g,x]$ is a homomorphism from $C_G(g/Z)$ to $Z$, whose kernel is $C_G(g)$. It follows that $C_G(g)$ is a subgroup of $C_G(g/Z)$ of index at most $|Z|$. The lemma follows.\ebew
\lem\label{index} Let $G$ be an \Mc-group. Then for any subgroup $H$ there is a bound for finite indices of the form $|H:C_H(g)|$.\elem
\bew Suppose not. Let $n$ be given by the \Mc-condition, and choose a maximal chain
$$H>C_H(g_0)>\cdots>C_H(g_0,\ldots,g_m)$$
with every group of finite index at least $n$ in its predecessor. Then $m<n$. However, if $\infty>|H:C_H(g)|>n\,|H:C_H(g_0,\ldots,g_m)|$, then 
$$|C_H(g_0,\ldots,g_m):C_H(g_0,\ldots,g_m,g)|>n,$$
contradicting maximality of $m$.\ebew
\defn Let $G$ be a group, and $H$, $K$ subgroups.\begin{itemize}
\item We say that $H$ is {\em almost contained} in $K$, written $H\ale K$, if $H\cap K$ has bounded index in $H$. Clearly, $\ale$ is transitive. If $H\ale K$ and $K\ale H$, then $H$ and $K$ are {\em commensurable}, denoted $H\aeq K$. 
\item The {\em almost centraliser} of $H$ in $K$ is the subgroup
$$\C_K(H)=\{g\in K: H\ale C_G(g)\}.$$
The {\em almost centre} $\Z(G)$ of $G$ is the characteristic subgroup $\Z(G)=\C_G(G)$.
\end{itemize}\edefn
\tats\label{f:mc}\begin{enumerate}\item If $K$ is definable and $H$ is type-definable, then by compactness $H\ale K$ if and only if $H\cap K$ has finite index in $H$.
\item \cite[Proposition 2.23]{He} In an \Mc-group the almost centraliser of a definable subgroup is definable by Lemma \ref{index}, and the almost centre is finite-by-abelian, as its conjugacy classes must be uniformly finite \cite{Neu}.\end{enumerate}\etats
\tats[{\cite[Theorem 2.10]{He}}]\label{symmetry} If $H$ and $K$ are type-definable, then $H\ale\C_G(K)$ if and only if $K\ale\C_G(H)$.\etats
\tats[{\cite[Theorem 2.18]{He}}]\label{commutator} Let $G$ be a group, $H$ and $K$ subgroups, and suppose
$$H\le N_G(K),\qquad H\le\C_G(K),\qquad\mbox{and}\qquad K\le\C_G(H)\mbox{ uniformly}$$
(meaning that there is $n<\omega$ such that $|H:C_H(k)|\le n$ for all $k\in K$). Then $[H,K]$ is finite.\etats
\kor Let $G$ be an \Mc-group, and $M$, $N$ normal subgroups of $G$. Then $[\C_M(N),\C_N(M)]$ is finite.\qed\ekor
\bem\label{central} It follows in particular that $F=[\C_G(\Z(G)),\Z(G)]$ is finite in an \Mc-group $G$, as $\Z(G)=\C_{\Z(G)}(G)$. But $G\ale\C_G(\Z(G))=\C_G(\C_G(G))$ by Fact \ref{symmetry}, and $F\le\Z(G)$, so $G_1=\C_G(\Z(G))\cap C_G(F)$ has finite index in $G$. Moreover, $F_1=G_1\cap F$ is finite central in $G_1$, and $(G_1\cap\Z(G))/F_1$ is central in $G_1/F_1$.\ebem

Recall that a group is {\em virtually} $P$ if it has a subgroup of finite index which is $P$.
\satz\label{abelian} Let $G$ be a broad pseudofinite dimensional \Mc-group satisfying strong fibration. Then $G$ has a $G$-broad definable finite-by-abelian subgroup $C$. More precisely, any $G$-broad minimal centralizer (up to finite index) of a finite tuple is virtually finite-by-abelian, and has a $G$-broad finite-by-abelian centraliser (of some bigger finite tuple) of finite index.\esatz
\bew By the \Mc-condition, there is a $G$-broad centralizer $C$ of some finite 
tuple, such that $C_C(g)$ is not $G$-broad for any $g\in C\setminus\Z(C)$.
Put $Z=\Z(C)$, a finite-by-abelian normal subgroup of $C$ which is definable by Fact \ref{f:mc}.

We claim that $Z$ is $G$-broad.
Otherwise $\dim_G(Z)=0$, and
$$\dim_G(C/Z)=\dim_G(C)-\dim_G(Z)=\dim_G(C)>0.$$
For $g\in C\setminus Z$ we have $\dim_G(C_C(g))=0$, whence for $\bar g=gZ$ we have
$$\begin{aligned}\dim_G(\bar g^{C/Z})&=\dim_G(g^CZ/Z)=\dim_G(g^CZ)-\dim_G(Z)\ge\dim_G(g^C)-\dim_G(Z)\\
&=\dim_G(C)-\dim_G(C_C(g))-\dim_G(Z)=\dim_G(C/Z).\end{aligned}$$
Hence for all $\bar g\in (C/Z)\setminus\{\bar1\}$ we have
$$\dim_G(C_{C/Z}(\bar g))=\dim_G(C/Z)-\dim_G(\bar g^{C/Z})=0.$$
As $C$ and $Z$ are definable, $C/Z$ is again 
pseudofinite, contradicting Proposition~\ref{centralizer}. This shows that
$\dim_G(Z)>0$, proving the claim. 

Now for $g\in \C_C(Z)$ the index $|Z:C_Z(g)|$ is finite. In particular $\dim(C_C(g))\ge\dim(Z)>0$ and $C_C(g)$ is $G$-broad. By minimality, $g\in\Z(C)$. It follows that $\Z(G)=\C_C(Z)$, which has finite index in $C$ since $Z\ale\C_C(C)$ implies $C\ale\C_C(Z)$ by Fact \ref{symmetry}. Moreover $\Z(C)$ is finite-by-abelian, and $C$ is virtually finite-by-abelian.

By compactness (or Lemma \ref{index}) there is a bound $n$ on the index $|C:C_C(c)|$ for $c\in\Z(C)$. For any $c'\in C\setminus\Z(C)$, the index
$|\Z(C):C_{\Z(C)}(c')|$ is infinite; let $c_0,\ldots,c_n$ in $\Z(C)$ lie in different cosets modulo $C_{\Z(C)}(c')$, and put $C'=C_C(c_i:i\le n)$, a subgroup of finite index in $C$. 

We claim that $|C:\Z(C)|>|C':\Z(C')|$. So suppose first that $c'\in C'\Z(C)$, say $c'=c''c$ with $c''\in C'$ and $c\in\Z(C)$. Then there are $i<j$ such that $c_i^{-1}c_j\in C_C(c)$, whence $c_i^{-1}c_j\in C_C(c''c)=C_C(c')$, a contradiction to the choice of the $c_i$.

It follows that $C'\Z(C)$ is a proper subgroup of finite index in $C$. Then $\Z(C'\Z(C))=\Z(C)$ and $\Z(C')=\Z(C)\cap C'$, whence 
$$|C:\Z(C)|>|C'\Z(C):\Z(C)|=|C':\Z(C')|,$$
proving the claim.

Inductively, we find a finite tuple $\bar c$ in $\Z(C)$ such that $C_C(\bar c)=\Z(C_C(\bar c))$, a $G$-broad finite-by-abelian centraliser of finite index in $C$.\ebew

Theorem \ref{abelian} holds in particular for any pseudofinite \Mc-group with coarse pseudofinite dimension (see Example \ref{example}(5)). Note that the \Mc-condition is just used in $G$, not in the section $C/Z$.

\kor\label{corTh} A superrosy pseudofinite group with $\Ut(G)\ge\omega^\alpha$ 
has a definable finite-by-abelian subgroup $A$ with 
$\Ut(A)\ge\omega^\alpha$.\ekor
\bew By \cite[Proposition 1.4]{EKP} a superrosy group is \Mc. If $\alpha$ is minimal with 
$\Ut(G)<\omega^{\alpha+1}$, put
$$\dim(X)\ge n\quad\text{if}\quad\Ut(X)\ge\omega^\alpha\cdot n,$$
then $\dim$ is an integer dimension with $0<\dim(G)<\infty$ and strong fibration (Example \ref{example}(2)). The 
assertion now follows from Theorem~\ref{abelian}.\ebew

\kor For any $d,d'<\omega$ there is $n=n(d,d')$ such that 
if $G$ is a finite group without elements $(g_i:i\le d')$ such that 
$$|C_G(g_i:i<j):C_G(g_i:i\le j)|\ge d$$ for all $j\le d'$, then $G$ has a 
subgroup $A$ with $|A'|\le n$ and $n\,|A|^n\ge|G|$.\ekor
\bew If the assertion were false, then given $d,d'$, there would be a sequence 
$(G_i:i<\omega)$ of finite groups satisfying the condition, such that $G_i$ has 
no subgroup $A_i$ with $|A_i'|\le i$ and $i\,|A_i|^i\ge|G_i|$. But any 
non-principal ultraproduct $G=\prod G_n/\U$ is a pseudofinite \Mc-group; by Proposition \ref{abelian} there is a definable subgroup $A$ with $A'$ finite and pseudofinite dimension
$\dim(A)\ge\frac1n\dim(G)$ for some $n<\omega$. Hence $\log|A_i|\ge\frac1n\log|G_i|-m$ for some $m\in\Z$ and almost all $i<\omega$, whence $e^{mn}|A_i|^n\ge|G|$.
For $i\ge\max\{n,|A'|,e^{mn}\}$ this yields a contradiction.\ebew

If $G$ is a definable group, a definable subgroup $H$ is {\em definably characteristic} in $G$ if it is invariant under all definable automorphisms of $G$.
\kor\label{dim1} Let $G$ be a pseudofinite \Mc-group of integer dimension $1$ with strong fibration. Then $G$ has a definably characteristic wide finite-by-abelian subgroup, which is a finite extension of the centraliser of a finite tuple (and hence quantifier-free definable).\ekor
\bew By Theorem \ref{abelian} there is a broad finite-by-abelian centraliser $C$ of a finite tuple. As $\dim(G)=1$ and $C$ is broad, $\dim(C)=1$ and $\dim(G/C)=0$. 

For any definable automorphism $\gamma$ of $G$ the image $C^\gamma$ is still definable as the centraliser of a finite tuple, and $\dim(G/C^\gamma)=0$. As $G/(C\cap C^\gamma)$ definably embeds into $G/C\times G/C^\gamma$, we have 
$$\dim(G/(C\cap C^\gamma))\le\dim(G/C)+\dim(G/C^\gamma)=0,$$
whence $\dim(C\cap C^\gamma)=1$. By Lemma \ref{index} there is a bound on the index of $C\cap C^\gamma$ in $C$ and in $C^\gamma$. Schlichting's Theorem now yields a definably characteristic subgroup $N$ commensurable with $C$, which is a finite extension of a finite intersection of conjugates of $C$ under definable automorphisms, and thus a finite extension of the centraliser of a finite tuple.

Now $\Z(N)$ is finite-by-abelian, and characteristic of finite index in $N$. It therefore contains all finite-by-abelian subgroups of finite index in $N$, and in particular $N\cap C$. So $\Z(N)$ is a finite extension of $N\cap C$, and thus definable as a finite extension of the centraliser of a finite tuple; clearly it is definably characteristic in~$G$.
\ebew

\section{Pseudofinite \Mc-groups of integer dimension $2$}
\satz\label{dim2} Let $G$ be a pseudofinite integer-dimensional \Mc-group with strong fibration. If $\dim(G)=2$, then $G$ has a broad definable
finite-by-abelian subgroup with wide normalizer.\esatz 
\bew Note that by Corollary \ref{dim1} we are done as soon as we find a definable subgroup of dimension $1$ with wide normalizer. By the \Mc-condition, there is a minimal wide centraliser of a finite tuple, up to finite index, and we can assume that this is already $G$. Thus 
$$\Z(G)=\{g\in G:\dim(C_G(c))=2\}.$$
If $\dim(\Z(G)\ge1$ we are done, taking $N=\Z(G)$. Otherwise $\Z(G)$ contains any subgroup $H$ of dimension $0$ with wide normaliser $N_G(H)$, since if $g\in H$ then $\dim(g^{N_G(H)})=0$ and $\dim(C_{N_G(H)}(g))=2$, whence $g\in\Z(G)$. By Remark \ref{central} there is a definable subgroup $G_1$ of finite index in $G$ and a finite central subgroup $F_1$ of $G_1$ such that $(G_1\cap\Z(G))/F_1$ is central in $G_1/F_1$. Now $G_1/F_1$ is still \Mc\ by Lemma \ref{byZ}; replacing $G$ by $G_1/F_1$ we can suppose $\Z(G)=Z(G)$.

We may assume $G=\prod_\U G_i$, where the $G_i$ are finite groups. For a subgroup $H\le G$ we shall put $\bar H=(HZ(G))/Z(G)$; similarly we shall put $\bar H_i=(H_iZ(G_i))/Z(G_i)$ for a subgroup $H_i\le G_i$. 
\beh There is some $k<\omega$ such that for almost all $i$ there is no direct product of simple nonabelian groups of length $k$ in $\bar G_i$.\ebeh
\bewbeh Suppose otherwise. Fix a decreasing chain $(I_k:k<\omega)$ of sets in $\U$ such that for $i\in I_k$ there is a direct product $\prod_{j<k}\bar N_{ij}\le\bar G_i$ of non-abelian simple groups. Fix such a product of length $k$ for all $i\in I_k\setminus I_{k+1}$, and note that $\bigcap_{k<\omega}I_k=\emptyset$, as $i\in I_k$ implies $k\le |G_i|$. Put $M_{ij}=N_{ij}'$. Then $M_{ij}Z(G_i)/Z(G_i)=\bar N_{ij}$, and $M_{ij}$ is perfect. Recall that the three subgroup Lemma states that for three normal subgroups $K,L,M$ of a group we have $[[K,L],M]\le[[L,M],K]\,[[M,K],L]$. Hence, for $j'\not=j$
$$[M_{ij},M_{ij'}]=[M_{ij}',M_{ij'}]\le[[M_{ij},M_{ij'}],M_{ij}]\,[[M_{ij'},M_{ij}],M_{ij}]=\{1\},$$
since $[M_{ij},M_{ij'}]\le Z(G_i)$.

For $j<\omega$ let $M_j=\prod_\U M_{ij}$; note that this is well-defined, and the $M_j$ are non-abelian and commute with one another. Put $H=\prod_{j<\omega}M_j$. Then if $m_j\in M_j\setminus Z(M_j)$ and $g_k=\prod_{j<k}m_j$, we have
$$|C_H(g_0,\ldots,g_\ell):C_H(g_0,\ldots,g_\ell,g_{\ell+n})|\ge 2^n$$
for all $\ell,n<\omega$, contradicting the \Mc-condition.
\ebew
\beh We may assume that almost all $\bar G_i$ are semisimple, i.e.\ have no non-trivial abelian normal subgroup.\ebeh
\bew Suppose otherwise, and let $\bar A_i$ be a non-trivial abelian normal subgroup of $\bar G_i$. Then for $a_i\in A_i\setminus Z(G_i)$ the element $a=[a_i]$ is not in $Z(G)$, so $\dim(C_G(a))\le 1$ and $\dim(a^G)\ge 1$. Consider $\bar A=Z(C_{\bar G}(\bar a^{\bar G}))$, an abelian normal subgroup of $\bar G$ containing $\bar a^{\bar G}$. Then its preimage $A$ is a definable normal subgroup of $G$ with $\dim(A)\ge\dim(a^G)\ge1$. Note that for $a\in A$ we have $a^{-1}a^A\subseteq Z(G)$ since $\bar A=A/Z(G)$ is abelian, so $\dim(a^A)=0$ and $\dim(C_A(a))=\dim(A)$. 

It follows that if $C$ is a minimal wide centraliser in $A$ of some finite tuple up to finite index (which exists by the \Mc-condition), then $C\aeq C^g$ for all $g\in G$, whence $\C_G(C)=\C_G(C^g)$ is normal in $G$. Therefore $\Z(\C_G(C))$ is a definable finite-by-abelian normal subgroup containing $C$, whence of dimension $\dim(A)\ge 1$: we are done.\ebew

Suppose $\bar G_i$ has at least two distinct minimal normal subgroups $\bar N_i$ and $\bar N'_i$ for almost all $i$. If $n_i\in N_i\setminus Z(G_i)$ for almost $i$, then $n=[n_i]\notin Z(G)$, so $\dim(n^G)\ge 1$. The same holds for $n'=[n'_i]$ with $n'_i\in N'_i\setminus Z(G_i)$. We put 
$$\bar N_0=C_{\bar G}(\bar n'^{\bar G})\supseteq \bar n^{\bar G}\qquad\mbox{and}\qquad
\bar N_1=C_{\bar G}(\bar N_0)\supseteq \bar n'^{\bar G}.$$
Then $\bar N_0$ and $\bar N_1$ are definable commuting normal subgroups of $\bar G$ of dimension at least $1$; their intersection must be trivial, as $\bar G$ is definably semisimple. Since $\dim(\bar G)=\dim(G)=2$ we must have $\dim(\bar N_j)=1$ for $j=0,1$. Now the preimages $N_j$ are definable normal subgroups of $G$ with $\dim(N_j)=\dim(\bar N_j)=1$; we are done again by Corollary \ref{dim1}.

So we may assume that almost each $\bar G_i$ has a unique minimal normal subgroup $\bar N_i$. Then $\bar N_i$ is a finite direct product of (finite) simple groups $\bar N_i^j$ for $j<k_i$, which are permuted transitively by $G_i$. By the first claim, almost all $k_i$ take the same value $k$.

Now every finite simple group is the product of two conjugacy classes \cite[Theorem 1.4]{GuMa}, so there are elements $\bar n_i,\bar n'_i$ with 
$$\bar N_i^0=\{\bar 1\}\cup\bar n_i^{\bar N_i^0}\,\bar n'^{\bar N_i^0}_i.$$
Note that if $\bar N_i^j=(\bar N_i^0)^{\bar g_{ij}}$ for some $\bar g_{ij}\in\bar G_i$, then
$$\bar N_i^j=\{\bar1\}\cup\bar n_i^{\bar N_i^0\bar g_{ij}}\,\bar n'^{\bar N_i^0\bar g_{ij}}_i.$$
Put $\bar X_i= \{\bar1\}\cup\bar n_i^{\bar G_i}\,\bar n'^{\bar G_i}_i$. It follows from normality that
$\bigcup_{j<k}\bar N_i^j\subseteq\bar X_i\subseteq\bar N_i$, so 
$\bar N_i^0= \{\bar1\}\cup\bar n_i^{\bar X_i}\,\bar n'^{\bar X_i}_i$ is uniformly definable, as is $\bar N_i^j$ for $j<k$.

For $j<k$ put $\bar N^j=\prod_\U \bar N_i^j$. Then $\bar N$ is definable as the direct product of the $\bar N^j$ for $j<k$. Since the $\bar N^j$ are definably isomorphic and $1\le \dim(\bar N)\le 2$, either $k=1$ and $\bar N=\bar N^1$, or $k=2=\dim(\bar N)$ and $\dim(\bar N^0)=\dim(\bar N^1)=1$. If $\dim(\bar N)=1$ we are done by Corollary \ref{dim1}; if $k=2$ the normalizer $N_G(\bar N^0)=N_G(\bar N^1)$ has index $2$ in $G$ and is wide, and we are done again.

So we may assume that $\dim(\bar N)=2$ and $\bar N$ is an ultraproduct of finite simple groups $\bar N_i$. As $\bar N$ is infinite, not almost all $\bar N_i$ can be sporadic, and we may assume they are all alternating or Chevalley groups (possibly twisted) over a finite field. But their rank (where the rank of the alternating group $A_k$ is $k$) must be bounded, as otherwise the $\bar N_i$ contain arbitrarily long direct products $\bar P_i$ of $A_4$ or PSL$_2$, contradicting the first claim.

It follows from \cite{Po99} that $\bar N$ must be a (possibly twisted) Chevalley group over a pseudofinite field. By results of Ryten \cite{Ry} (see also \cite{EM}) the pure group $\bar N$ is bi-interpretable with a pseudofinite (difference) field $F$ of $SU$-rank $1$. Now $\bar N$ is internal in $F$, and $F$ is internal in any generic definable subset $X$ of $F$, whence 
$$\dim(X)\ge\frac1m\dim(F)\ge\frac1n\dim(\bar N)$$
for some positive integers $m,n$,
and $\dim(X)>0$. As the dimension is integer, this means $\dim(X)\ge 1$ for any definable generic subset of $F$. But $SU(\bar N)\ge 3$, so a generic element of $\bar N$ (interpreted in $F$) contains at least three independent generic coordinates from $F$; as $SU$-rank is definable in $F$, we obtain a definable subset of $\bar N$ of dimension $3$, contradicting $\dim(\bar N)\le\dim(G)=2$.\ebew

\kor\label{soluble} Let $G$ be a pseudofinite group whose definable sections 
are \Mc, and $\dim$ an integer dimension on $G$ with strong fibration. If 
$\dim(G)=2$, then $G$ has a definable wide soluble subgroup.\ekor
\bew By Theorem \ref{dim2}, there is a definable finite-by-abelian group $N$ 
such that $N_G(N)$ is wide. Replacing $N$ by $C_N(N')$, we may assume that 
$N$ is (finite central)-by-abelian. If $\dim(N)=2$ we are done. Otherwise 
$\dim(N_G(N)/N)=1$; by 
Corollary \ref{dim1} there is a definable
finite-by-abelian subgroup $S/N$ with $\dim(S/N)=1$. As above 
we may assume that $S/N$ is (finite central)-by-abelian, so $S$ is soluble. 
Moreover, 
$$\dim(S)=\dim(N)+\dim(S/N)=1+1=2,$$
so $S$ is wide in $G$.\ebew

\kor A pseudofinite superrosy group $G$ with 
$\omega^\alpha\cdot2\le\Ut(G)<\omega^\alpha\cdot3$ has a definable soluble 
subgroup $S$ with $\Ut(S)\ge\omega^\alpha\cdot2$.\ekor
\bew Superrosiness implies that all definable sections of $G$ are \Mc. We put
$$\dim(X)=n\quad\Leftrightarrow\quad\omega^\alpha\cdot 
n\le\Ut(X)<\omega^\alpha\cdot(n+1).$$
This defines an integer dimension with strong fibration, and $\dim(G)=2$. The result now follows 
from Corollary \ref{soluble}.\ebew

\end{document}